\documentclass[11pt]{amsart}
\usepackage{fullpage}
\usepackage{color}
\usepackage{graphicx,psfrag} 
\usepackage{color} 
\usepackage{tikz}
\usepackage{pgffor}
\usepackage{todonotes}

\usepackage{perpage}     

\usepackage{wrapfig}
\usepackage{enumerate}
\usepackage{subfigure}
\usepackage{caption}

\usepackage[normalem]{ulem}
\usepackage[makeroom]{cancel} 

\newtheorem{theorem}{Theorem}[section]

\newtheorem{proposition}[theorem]{Proposition}

\theoremstyle{definition}

\newtheorem{example}[theorem]{Example}
\theoremstyle{remark}

\newtheorem{question}[theorem]{Question}

\newtheorem{ex-prop}[theorem]{Example-Proposition}

\numberwithin{equation}{section}

\definecolor{gray}{rgb}{.5,.5,.5}

\definecolor{black}{rgb}{0,0,0}

\definecolor{blue}{rgb}{0,0,1}

\definecolor{red}{rgb}{1,0,0}
\def\red{\color{red}}
\definecolor{green}{rgb}{0,1,0}

\definecolor{yellow}{rgb}{1,1,.4}

\definecolor{purple}{rgb}{1,0,1}

\definecolor{gold}{rgb}{.5,.5,.2}

\definecolor{darkgreen}{rgb}{0,.5,0}

\definecolor{greenbean}{RGB}{199, 237, 204}

\definecolor{RED}{rgb}{1,0,0}

\newcommand{\G}{\mathbb{G}}

\begin{document}

\title{A counterexample to Question 1 of\\ ``A survey on the Turaev genus of knots''}

\author{Cody Armond}
\address{Department of Mathematics\\
University of Iowa\\
Iowa City, Iowa, USA}
\email{cody-armond@uiowa.edu}

\author{Moshe Cohen}
\address{Department of Mathematics, Technion -- Israel Institute of Technology, Haifa 32000, Israel}
\email{mcohen@tx.technion.ac.il}

\date{\today}

\begin{abstract}
In ``A survey on the Turaev genus of knots,'' Champanerkar and Kofman propose several open questions.  The first asks whether the polynomial whose coefficients count the number of quasi-trees of the all-A ribbon graph obtained from a diagram with minimal Turaev genus is an invariant of the knot.   We answer negatively by showing a counterexample obtained from the two diagrams of $8_{21}$ on the KnotAtlas and KnotScape.
\end{abstract}


\maketitle

\section{Introduction}
\label{sec:Intro}

Champanerkar and Kofman offer a very complete ``survey on the Turaev genus of knots'' \cite{ChKo:dessin}, and we defer the reader to this short survey rather than repeat most of the background here.

In an earlier work with Stoltzfus, they define a polynomial whose coefficients count the number of quasi-trees of the all-A ribbon graph $\G$ obtained from a diagram with minimal Turaev genus.  This comes from the Bollob\'{a}s-Riordan-Tutte polynomial $C(\G,X,Y,Z)$.

\begin{proposition}
\label{prop:q}
\cite[Proposition 3.2]{ChKoSt} Let $q(\G;t, Y ) = C(\G; 1, Y, tY^{-2})$. Then $q(\G;t, Y )$ is a polynomial in $t$ and $Y$ such that
\begin{equation}
\label{eq:q}
q(\G;t) := q(\G;t, 0) = \sum_j a_j t^j
\end{equation}
where $a_j$ is the number of quasi-trees of genus $j$. Consequently, $q(\G; 1)$ equals the number of quasi-trees of $\G$.
\end{proposition}

Dasbach, Futer, Kalfagianni, Lin, and Stoltzfus in \cite[Theorem 3.2]{DaFuKaLiSt2} show that the evaluation of this polynomial with $t=-1$ gives the determinant of the knot.

The recent survey paper asks whether the polynomial itself is an invariant when the Turaev genus $g_T$ of the diagram is equal to that of the knot, that is, when it is minimal.

\begin{question}
\label{question}
\cite[Question 1]{ChKo:dessin} 
Let $\G$ be the all-A ribbon graph for a diagram $D$ of a knot $K$. If $g_T(D) = g_T(K)$, is $q(\G;t)$ an invariant of $K$?
\end{question}

We give a negative answer to this question by providing a counterexample.

\begin{theorem}
\label{thm:main}
The polynomial whose coefficients count the number of quasi-trees of the all-A ribbon graph obtained from diagram with minimal Turaev genus is not an invariant of the knot.
\end{theorem}

We prove this Theorem \ref{thm:main} by considering the two diagrams of $8_{21}$ obtained from the Knot Atlas \cite{knotatlas} and KnotScape \cite{knotscape}.  We address these cases in Examples \ref{ex:atlas} and \ref{ex:info}, respectively.

We rely on an algorithm given by Armond, Druivenga, and Kindred in \cite{ArDrKi} to obtain alternating diagrams on a surface with minimal Turaev genus.

\section{A counterexample:  diagrams from the Knot Atlas and KnotScape}
\label{sec:Counterexample}

\begin{figure}[h]
\centering
\subfigure{
    \includegraphics[scale=0.4]{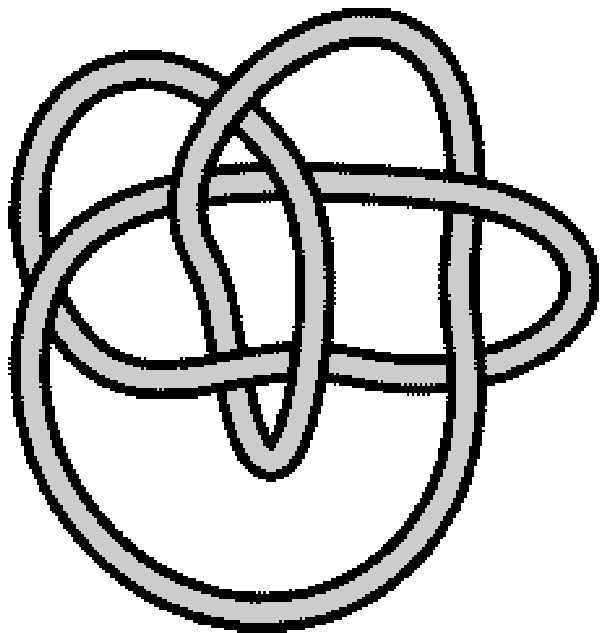}
}
\hspace{2cm}
\subfigure{
    \includegraphics[width=.3\textwidth, height=3cm]{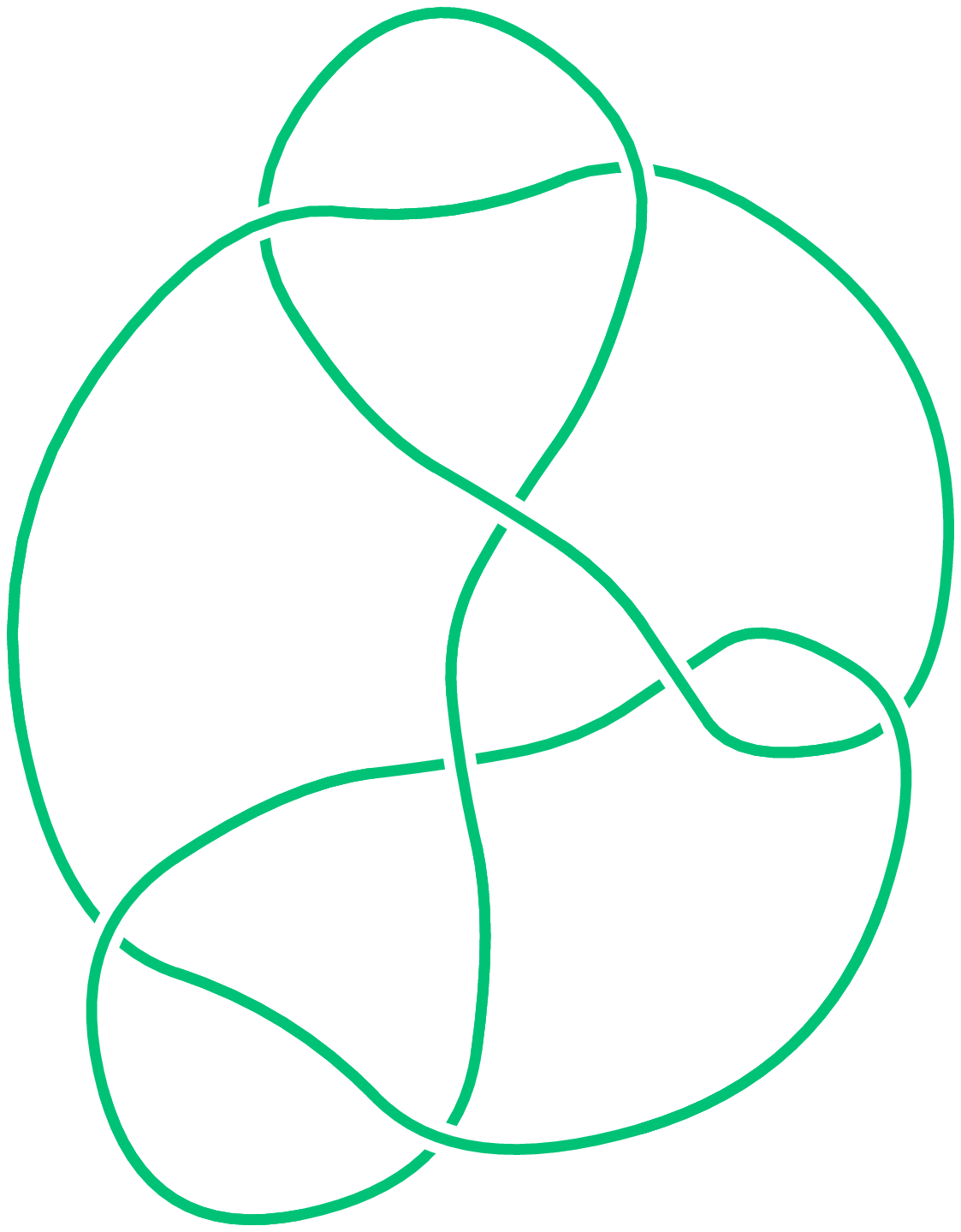}
}
\caption{\label{fig:821} The knot $8_{21}$ presented in diagrams given by the Knot Atlas \cite{knotatlas} and KnotScape \cite{knotscape}.}
\end{figure}

In Examples \ref{ex:atlas} and \ref{ex:info} below, we count the quasi-trees of the all-A ribbon graph obtained from diagrams coming from the Knot Atlas \cite{knotatlas} and KnotScape \cite{knotscape}, respectively, as shown in Figure \ref{fig:821}.  We show that the polynomial $q(\G,t)$ is not invariant on the knot.


\begin{example}
\label{ex:atlas}
Consider first the knot diagram of $8_{21}$ given by the Knot Atlas \cite{knotatlas}, as shown in Figure \ref{fig:821}.  This diagram has Turaev genus 2.  We perform a Reidemeister III move on the upper central three crossings to obtain a diagram of Turaev genus 1.
\end{example}

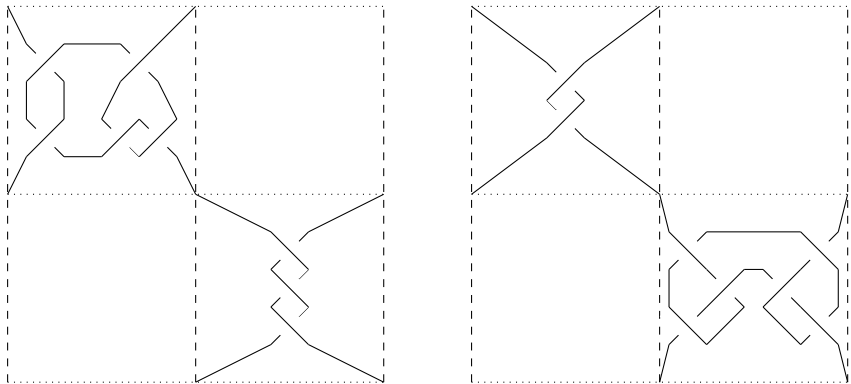
\begin{figure}[h]
\centering
\subfigure{
\begin{tikzpicture}[scale=.5]
\draw[dashed] (0,0) -- (0,10);
\draw[dashed] (5,0) -- (5,10);
\draw[dashed] (10,0) -- (10,10);
\draw[dotted] (0,0) -- (10,0);
\draw[dotted] (0,5) -- (10,5);
\draw[dotted] (0,10) -- (10,10);

\foreach \x/ \y in {.5/6, .5/8, 2.5/6, 3.5/6, 3/8}
    {
    \draw (\x+1,\y) -- (\x,\y+1);
    \draw[color=white, line width=10] (\x,\y) -- (\x+1,\y+1);
    \draw (\x,\y) -- (\x+1,\y+1);
    }

\foreach \x/ \y in {7/1, 7/2, 7/3}
    {
    \draw (\x,\y) -- (\x+1,\y+1);
    \draw[color=white, line width=10] (\x+1,\y) -- (\x,\y+1);
    \draw (\x+1,\y) -- (\x,\y+1);
    }

\draw (.5,7) -- (.5,8);
\draw (1.5,7) -- (1.5,8);
\draw (1.5,6) -- (2.5,6);
\draw (1.5,9) -- (3,9);
\draw (2.5,7) -- (3,8);
\draw (4.5,7) -- (4,8);

\draw (5,0) -- (7,1);
\draw (10,0) -- (8,1);
\draw (5,5) -- (7,4);
\draw (10,5) -- (8,4);
\draw (0,5) -- (.5,6);
\draw (5,5) -- (4.5,6);
\draw (0,10) -- (.5,9);
\draw (5,10) -- (4,9);
\end{tikzpicture}
}
\hspace{2em}
\subfigure{
\begin{tikzpicture}[scale=.5]
\draw[dashed] (0,0) -- (0,10);
\draw[dashed] (5,0) -- (5,10);
\draw[dashed] (10,0) -- (10,10);
\draw[dotted] (0,0) -- (10,0);
\draw[dotted] (0,5) -- (10,5);
\draw[dotted] (0,10) -- (10,10);

\foreach \x/ \y in {2/6.5, 2/7.5, 6.25/2, 7.75/2}
    {
    \draw (\x+1,\y) -- (\x,\y+1);
    \draw[color=white, line width=10] (\x,\y) -- (\x+1,\y+1);
    \draw (\x,\y) -- (\x+1,\y+1);
    }

\foreach \x/ \y in {5.25/1, 5.25/3, 8.75/1, 8.75/3}
    {
    \draw (\x,\y) -- (\x+1,\y+1);
    \draw[color=white, line width=10] (\x+1,\y) -- (\x,\y+1);
    \draw (\x+1,\y) -- (\x,\y+1);
    }

\draw (6.25,1) -- (7.25,2);
\draw (7.75,2) -- (8.75,1);
\draw (5.25,2) -- (5.25,3);
\draw (9.75,2) -- (9.75,3);
\draw (7.25,3) -- (7.75,3);
\draw (6.25,4) -- (8.75,4);

\draw (5,0) -- (5.25,1);
\draw (10,0) -- (9.75,1);
\draw (5,5) -- (5.25,4);
\draw (10,5) -- (9.75,4);
\draw (0,5) -- (2,6.5);
\draw (5,5) -- (3,6.5);
\draw (0,10) -- (2,8.5);
\draw (5,10) -- (3,8.5);
\end{tikzpicture}
}
\caption{\label{fig:Window} Alternating diagrams on the torus for $8_{21}$ coming from the KnotAtlas and KnotScape, respectively, after applying the algorithm of \cite{ArDrKi}.}
\end{figure}

Armond, Druivenga, and Kindred \cite{ArDrKi} give an algorithm to obtain an alternating diagram on a surface.  We apply this to obtain a Heegaard diagram, where the dashed and dotted lines represent $\alpha$ and $\beta$ curves, respectively, as given on the left-hand side in Figure \ref{fig:Window}.

\begin{figure}[h]
\centering
\subfigure{
\begin{tikzpicture}

\draw (2,3) ellipse (2.1cm and .6cm);
\draw (2,3) ellipse (2cm and .5cm);
\draw (4,2) ellipse (.5cm and 2cm);

	\fill[color=black] (3.55,2.66) circle (5pt);
	\fill[color=black] (3.49,1.69) circle (5pt);
	\fill[color=black] (3.65,.69) circle (5pt);
\draw (3.85,2.19) node {$c$};
\draw (3.85,1.19) node {$b$};
\draw (4,.4) node {$a$};

	\fill[color=black] (.45,2.66) circle (5pt);
\draw (2,3.85) node {$h$};
\draw (2,3.25) node {$g$};

	\fill[color=black] (2,2.3) circle (5pt);
\draw (1.25,2.15) node {$e$};
\draw (2.75,2.15) node {$d$};
\draw (2,2.75) node {$f$};
\end{tikzpicture}
}
\hspace{2em}
\subfigure{
    \includegraphics[width=.3\textwidth, height=4cm]{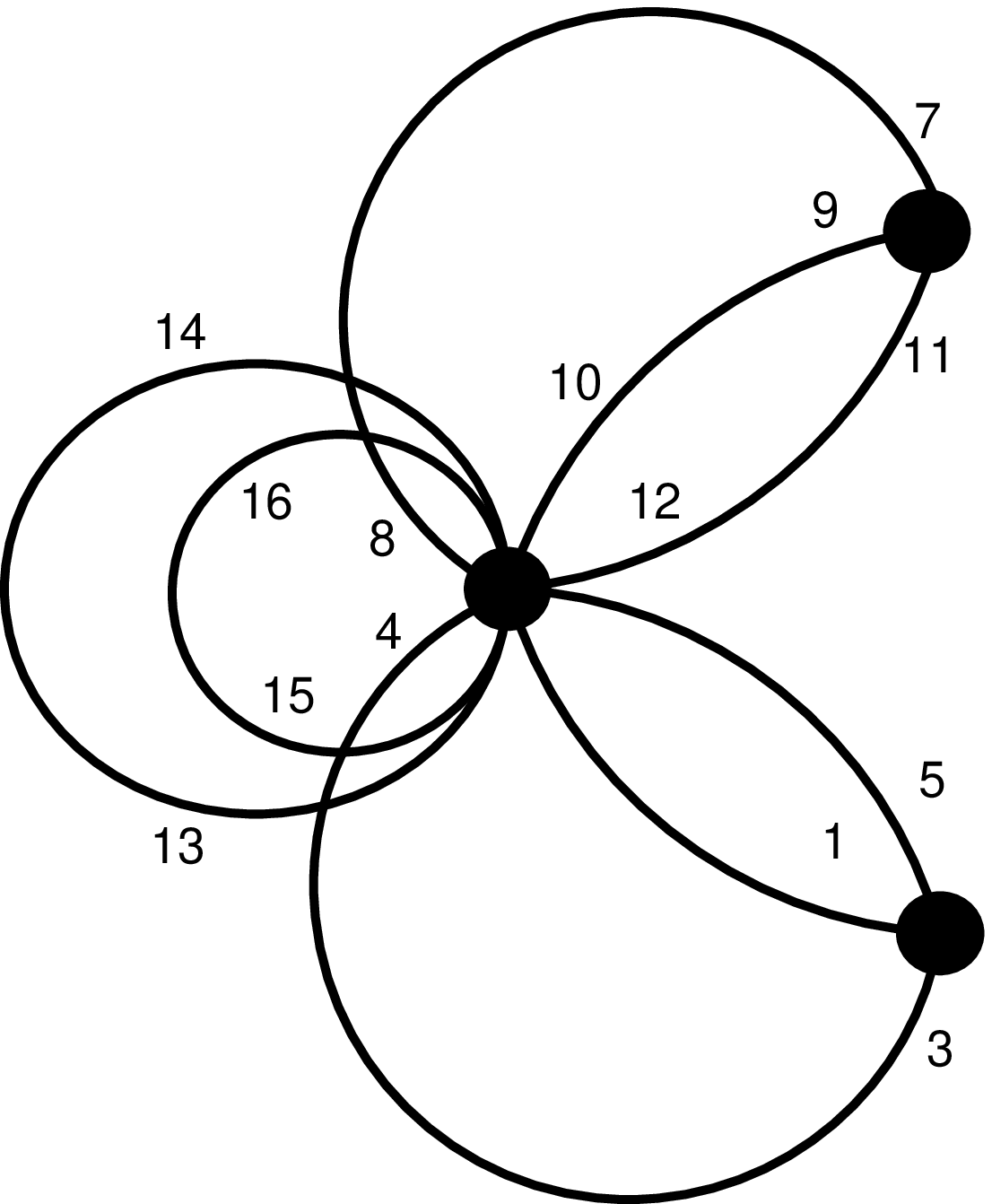}
}
\caption{\label{fig:Dessin} The all-A ribbon graphs for diagrams of $8_{21}$ coming from the KnotAtlas after a Reidemeister move III and from KnotScape (as appearing in \cite[Figure 3]{DaFuKaLiSt2}), respectively.}
\end{figure}

We checkerboard color this diagram on the torus to obtain the all-A ribbon graph given on the left-hand side in Figure \ref{fig:Dessin}.  We proceed to count the number of quasi-trees.

First of all, any spanning tree of $\G$ must contain exactly two edges from the loop consisting of edges $a$, $b$, and $c$.  From here the spanning trees fall into two classes: those with one of the two edges $g$ and $h$ and those with neither $g$ nor $h$. Any spanning tree in the first class must contain one of the two edges $d$ and $e$ giving a total of $3 \times 2 \times 2 = 12$ spanning trees in the first class. Any spanning tree in the second class must contain two of the three edges $d$, $e$, and $f$ giving a total of $3 \times 3 = 9$ spanning trees in the second class. Thus for this ribbon graph we get $a_0=21$.

A quasi-tree of $\G$ with genus $1$ must contain all of the edges $a$, $b$, and $c$ as well as one of the two edges $g$ and $h$ and again these quasi-trees fall into two classes: those that contain the edges $d$ and $e$ but not $f$ and those that contain the edge $f$ and exactly one of the edges $d$ and $e$. This gives us $a_1 = 4+2=6$ for this ribbon graph.

Thus, we obtain $q(\G,t)=6t+21$.

\begin{example}
\label{ex:info}
Now consider the knot diagram of $8_{21}$ given by KnotScape \cite{knotscape} having Turaev genus 1 already and appearing on the right-hand side of Figure \ref{fig:821}.
\end{example}

We apply the algorithm of \cite{ArDrKi} to obtain an alternating diagram on a surface, which again is a Heegaard diagram, where the dashed and dotted lines represent $\alpha$ and $\beta$ curves, respectively, as given on the right-hand side in Figure \ref{fig:Window}.


We checkerboard color this diagram on the torus to obtain 
 the all-A ribbon graph, given in Figure 3 from \cite{DaFuKaLiSt2}, which we include on the right-hand side in our Figure \ref{fig:Dessin}.

As observed in \cite{DaFuKaLiSt2}, this ribbon graph contains 9 spanning trees and 24 genus-1 quasi-trees, yielding $q(\G,t)=24t+9$.


\medskip

\textbf{Acknowledgements.} The second author would love to thank the first author and his wife for hosting him in Iowa City while this work was completed.


\newcommand{\etalchar}[1]{$^{#1}$}
\def\cprime{$'$}
\providecommand{\bysame}{\leavevmode\hbox to3em{\hrulefill}\thinspace}
\providecommand{\MR}{\relax\ifhmode\unskip\space\fi MR }
\providecommand{\MRhref}[2]{%
  \href{http://www.ams.org/mathscinet-getitem?mr=#1}{#2}
}
\providecommand{\href}[2]{#2}

\end{document}